\documentclass[11pt,leqno]{article}
\usepackage{latexsym}
\usepackage{amsmath,amssymb}
\usepackage{amsthm}
\usepackage[spanish,english]{babel}
\title{Extension of Carter subgroups in $\pi$-separable groups}
\author{M. Arroyo-Jord\'a, P. Arroyo-Jord\'a, R. Dark,\\ A.~D. Feldman, and M.~D. P\'{e}rez-Ramos}

\date{}

\newtheorem{teo}{Theorem}[section]
\newtheorem{lem}[teo]{Lemma}
\newtheorem{pro}[teo]{Proposition}
\newtheorem{cor}[teo]{Corollary}
\theoremstyle{definition}
\newtheorem{rem}[teo]{Remark}

\newtheorem{de}[teo]{Definition}

\DeclareMathOperator{\Syl}{Syl}
\DeclareMathOperator{\Hall}{Hall}

\DeclareMathOperator{\Char}{Char}
\DeclareMathOperator{\Alt}{Alt}

\DeclareMathOperator{\proj}{Proj}
\DeclareMathOperator{\cov}{Cov}
\DeclareMathOperator{\p}{\mathbb{P}}

\begin{document}
\maketitle

\begin{abstract}

Let $\pi$ be a set of primes. We show that $\pi$-separable groups have a conjugacy class of subgroups which specialize to Carter subgroups, i.e. self-normalizing nilpotent  subgroups, or equivalently, nilpotent projectors, when specializing to soluble groups.
\medskip

\medskip

\end{abstract}

\section{Introduction: From soluble  to $\pi$-separable groups}

All groups considered are finite.

We pursue an extension of the theory of finite soluble groups to the universe of $\pi$-separable groups, $\pi$ a set of primes. Classical results of Hall theory state that soluble groups are characterized by the existence of Hall $\rho$-subgroups for all sets of primes $\rho$ (P. Hall~\cite{H28,H37}). When the set of primes $\pi$ is fixed, $\pi$-separable groups have Hall $\pi$-subgroups, and also every $\pi$-subgroup is contained in a conjugate of any Hall $\pi$-subgroup (S. A. \v{C}unihin~\cite{Cu}). The main results in this paper (Section~\ref{sCarter-like}, Theorems~\ref{teo1}~and \ref{teo2}) show that $\pi$-separable groups have a conjugacy class of subgroups which specialize to Carter subgroups, i.e. self-normalizing nilpotent  subgroups, or equivalently, nilpotent projectors, when specializing to soluble groups. In the theory of soluble groups, Carter subgroups are the cornerstone for the existence and conjugacy of injectors associated to Fitting classes. In a forthcoming paper \cite{ADFP}, our Carter-like subgroups are used to generalize these results to  $\pi$-separable groups.

With this aim we analyze first the reach of $\pi$-separability further from the universe of soluble groups. We refer to \cite{Go} for basic results on $\pi$-separable groups, and to \cite{DH} for background on classes of groups;  we shall adhere to their notations. If $\pi$ is a set  of primes, let us recall that a group $G$ is $\pi$-separable if every composition factor of $G$ is either a $\pi$-group or a $\pi'$-group, where $\pi'$ stands for the complement of $\pi$ in the set $\p$ of all prime numbers. We notice that a group $G$ is soluble if and only if it is $\rho$-separable for all sets of primes $\rho$. Regarding $\pi$-separability, $\pi$ a set of primes, it is clearly equivalent to $\pi'$-separability, so that there is no loss of generality to assume that $2\in \pi$. Then, by the Feit-Thompson theorem, a $\pi$-separable group is $\pi'$-soluble, i.e. the group is $\pi'$-separable with every $\pi'$-composition factor a $p$-group for some prime $p\in \pi'$. The following extension  for $\pi$-separable groups is  easily proved:

\begin{pro}\label{pro0} For a group $G$, if $2\in \pi\subseteq \p$ the following statements are pairwise equivalent:
\begin{enumerate} \item $G$ is $\pi$-separable;
\item $G$ is $\rho$-separable for every set of primes $\rho$ such that either $\pi \subseteq \rho$ or $\pi\cap \rho =\emptyset$;
\item $G$ is $\pi'$-separable ($\pi'$-soluble).
\end{enumerate}
\end{pro}

\begin{rem} The need of the hypothesis $2\in \pi$ for the validity of the equivalences in Proposition~\ref{pro0} is clear. Certainly, for a group $G$ and any set of primes $\pi$, it holds that  $2\rightarrow 1 \leftrightarrow 3$. But statement 2 implies that the $\pi'$-compositions factors of the group $G$ are soluble. Then, by the Feit-Thompson theorem, if $2\notin \pi$, statement 2 is equivalent to the solubility of the group $G$ and $1\nrightarrow 2$ in general.
\end{rem}

Consequently, by Proposition~\ref{pro0}, if $2\in \pi\subseteq \p$, every $\pi$-separable group possesses a Sylow $p$-complement of $G$, i.e. a Hall $p'$-subgroup of $G$, for each $p\in \pi'$, as well as a Hall $\pi'$-subgroup, and these  are pairwise permutable subgroups with coprime indices in the group. In analogy with \cite[I. Definitions (4.1), (4.5), (4.7)]{DH} it appears to be natural now to introduce the following concepts, which are proven to hold  in  $\pi$-separable groups if $2\in \pi$, by the previous comment and as explained below:

\begin{de} Let $G$ be a group and $\pi$ be a set of primes.
\begin{itemize}\item[$\mathbf{K}_\pi$:] A \emph{complement $\pi$-basis} of  $G$ is a set $\mathbf{K}_\pi$ containing exactly one Sylow $p$-complement of $G$, i.e. a Hall $p'$-subgroup of $G$, for each $p\in \pi'$, and exactly one Hall $\pi'$-subgroup.
\item[$\mathbf{\Sigma_\pi}$:] A \emph{Hall $\pi$-system} of  $G$ is a set $\mathbf{\Sigma_\pi}$ of Hall subgroups of $G$ satisfying the following two properties:
    \begin{enumerate}
    \item For each set of primes $\rho$ such that either $\pi \subseteq \rho$ or $\pi\cap \rho =\emptyset$, $\mathbf{\Sigma_\pi}$ contains exactly one Hall $\rho$-subgroup.
    \item If $H,K\in \mathbf{\Sigma_\pi}$, then $HK=KH$.
    \end{enumerate}
    \item[$\mathbf{B_\pi}$:] A \emph{Sylow $\pi$-basis} of  $G$ is a set $\mathbf{B_\pi}$ of subgroups of $G$ satisfying the following two properties:
    \begin{enumerate}
    \item $\mathbf{B_\pi}$  contains exactly one Hall $\pi$-subgroup, and exactly one Sylow $p$-subgroup for each $p\in \pi'$.
    \item If $H,K\in \mathbf{B_\pi}$, then $HK=KH$.
    \end{enumerate}
\end{itemize}
\end{de}

Obviously these systems may not exist in arbitrary  groups. By the previous comments, if the group is $\pi$-separable and $2\in \pi$, then complement $\pi$-bases do exist.

Note also that Sylow $\pi$-systems and  complement and Sylow  $\pi$-bases are hereditary with respect to normal subgroups and factor groups.

For any set $\rho$ of primes and a group $G$, we denote by $\Hall_{\rho}(G)$ the set of all Hall $\rho$-subgroups of $G$. If $p$ is a prime, then $\Syl_p(G)$ will denote the set of all Sylow $p$-subgroups of $G$. We keep mimicking the exposition in \cite[Section I.4]{DH}. The arguments there are easily adapted to prove the following corresponding results.

\begin{pro}\label{proKH}\textup{\cite[I. Proposition (4.4)]{DH}} Assume that the group $G$ has a complement $\pi$-basis, say $\mathbf{K_\pi}$ (particularly, if the group $G$ is  $\pi$-separable  and  $2\in \pi\subseteq \p$).  If $\rho$ is a set of primes such that $\pi\subseteq \rho$, let $G_\rho=\bigcap\{X\mid X\in \Hall_{p'}(G)\cap \mathbf{K_\pi},\ p\in \rho'\subseteq \pi'\}$. On the other hand, if $\rho$ is a set of primes such that $\pi\cap \rho=\emptyset$, let
$G_\rho=\bigcap\{X\mid (X\in \Hall_{p'}(G)\cap \mathbf{K_\pi},\ p\in \rho'\cap \pi')\ \vee\ (X\in \Hall_{\pi'}(G)\cap \mathbf{K_\pi})\}$. Then
\begin{enumerate}
\item $\mathbf{\Sigma_\pi}:= \{G_\rho\mid (\pi\subseteq \rho\subseteq \p)\ \vee (\rho \subseteq \p,\ \rho \cap \pi=\emptyset)\}$ is a Hall $\pi$-system of $G$, and
    \item $\mathbf{\Sigma_\pi}$ is the unique Hall $\pi$-system of $G$ containing $\mathbf{K_\pi}$.
\end{enumerate}

We shall say that $\mathbf{\Sigma_\pi}$ is the Hall $\pi$-system \emph{generated} by the complement $\pi$-basis $\mathbf{K_\pi}$.
\end{pro}

\begin{cor}\textup{\cite[I. Corollary (4.6)]{DH}} Let $G$ be a  $\pi$-separable group,  $2\in \pi\subseteq \p$. Then there is a bijective map between the set of all complement $\pi$-bases and the set of all Hall $\pi$-systems of $G$, such that to each complement $\pi$-basis corresponds the Hall $\pi$-system generated by it, and conversely, to each Hall $\pi$-system corresponds the complement $\pi$-basis contained in it.
\end{cor}

On the other hand, it is clear that every Hall $\pi$-system contains a unique Sylow $\pi$-basis. Also, each Sylow $\pi$-basis \emph{generates} a unique Hall $\pi$-system, by taking the product of the suitable elements in the basis to construct each element in the Hall $\pi$-system. We can easily state also the following result:

\begin{cor}\textup{\cite[I. Lemma (4.8)]{DH}} Let $G$ be a  $\pi$-separable group,  $2\in \pi\subseteq \p$. Then
there is a bijective map between the set of all Hall $\pi$-systems and the set of all Sylow $\pi$-bases of $G$, such that to each Hall $\pi$-system corresponds the Sylow $\pi$-basis contained in it, and conversely, to each Sylow $\pi$-basis corresponds the Hall $\pi$-system generated by it as described above.
\end{cor}

If $\mathbf{K_\pi}$ is a complement $\pi$-basis of a group $G$ and $g\in G$, it is clear that $\mathbf{K_\pi}^g:=\{X^g\mid X\in \mathbf{K_\pi}\}$ is again a complement $\pi$-basis of $G$, and this defines an action by conjugation of $G$ on the set of all complement $\pi$-bases of $G$. Analogously, any group $G$ acts by conjugation on the  set of all its Hall $\pi$-systems, and also on the set of all its Sylow $\pi$-bases.
\smallskip

By Proposition~\ref{pro0} again we have that if $G$ is a $\pi$-separable group, $2\in \pi\subseteq \p$, then $G$ acts transitively on the set of all Hall $p'$-subgroups for every $p\in \pi'$, as well as on the set of all Hall $\pi'$-subgroups, and the following result also holds.

\begin{teo}\textup{\cite[I. Theorems (4.9), (4.10), (4.11), Corollary (4.12)]{DH}} Let $G$ be a $\pi$-separable group, $2\in \pi$. Then:
 \begin{enumerate} \item The number of Hall $\pi$-systems of $G$ is $\prod_{S\in \mathbf{K_\pi}}|G:N_G(S)|$, where $\mathbf{K_\pi}$ is a complement $\pi$-basis of $G$.
 \item The group $G$ acts transitively by conjugation on the set of all complement $\pi$-bases, on the set of all Hall $\pi$-systems, as well as on the set of all Sylow $\pi$-bases.
     \end{enumerate}
\end{teo}

One might wish that the existence of Hall $\pi$-systems in finite groups would characterize $\pi$-separability, but this is not the case, even assuming the transitive action of the group by conjugation on the set of all Hall $\pi$-systems. The alternating group of degree $5$ together with the set $\pi=\{2,3\}$ is a counterexample. A characterization of $\pi$-separability by the existence of Hall subgroups had been in fact given by Z. Du~\cite{Du}, as shown in the next result. For notation, for a group $G$ and any set of primes $\rho$, the group $G$ is said to satisfy $E_{\rho}$ if $G$ has a Hall $\rho$-subgroup; if in addition each $\rho$-subgroup is contained in the conjugate of a Hall $\rho$-subgroup, it is said that $G$ satisfies $D_\rho$.

\begin{teo}\textup{\cite[Theorems 1, 3]{Du}} For a group $G$ and a set of primes $\pi$, the following statements are pairwise equivalent:
\begin{itemize}\item[(i)] $G$ is $\pi$-separable.
\item[(ii)] $G$ satisfies:
\begin{enumerate}
\item $E_\pi$ and $E_{\pi'}$;
\item $E_{\pi\cup \{q\}}$ and $E_{\pi'\cup \{p\}}$, for all $p\in \pi$, $q\in \pi'$.
\end{enumerate}
\item[(iii)] $G$ satisfies:
\begin{enumerate}
\item $E_\pi$ and $E_{\pi'}$;
\item $E_{\{p,q\}}$, for all $p\in \pi$, $q\in \pi'$.
\end{enumerate}
\end{itemize}
\end{teo}

We point out finally that the existence of Hall $\pi$-systems together with $\pi$-dominance do characterize $\pi$-separability, as we prove next.

\begin{rem}\label{dpidpi'} We notice that for any group $G$ and any $\pi\subseteq \p$, the existence of complement $\pi$-bases is equivalent to the existence of Hall $\pi$-systems, and also to the existence of Sylow $\pi$-bases, by Proposition~\ref{proKH} and the corresponding constructions. In this case the group $G$ satisfies $E_\pi$ and $E_{\pi'}$, and then, if $2\in \pi$, $G$ satisfies $D_{\pi'}$ (see \cite{AF}).
\end{rem}

For any positive integer $n$, we denote by $\pi(n)$ the set of primes dividing $n$; for  the order $|G|$ of a
group G, we set  $\pi(G)=\pi(|G|)$.

\begin{teo} Assume that the group $G$ has a complement $\pi$-basis and satisfies $D_\pi$, where $2\in \pi\subseteq \p$. Then $G$ is $\pi$-separable.
\end{teo}
\smallskip

{\it Proof.} We argue by induction on the order of $G$. Let $N$ be a normal subgroup
 of $G$. By \textup{(\cite[Theorem 7.7]{RV})}, $G$ satisfies $D_\pi$ if and only if $N$ and $G/N$ satisfies $D_\pi$. We may then assume that $G$ is a simple group. If $|\pi(G)\cap \pi'|\ge 2$,  then $G$ would satisfy $E_{p'}$ and $E_{q'}$, with $p,q$ odd different primes dividing the order of $G$, which would imply that $G$ would not be simple by \cite[Corollary 5.5]{AF}. Consequently we may assume that $|\pi(G)\cap \pi'|= 1$. On the other hand, by hypothesis and Remark~\ref{dpidpi'}, $G$ satisfies $D_\pi$ and $D_{\pi'}$, and we may assume that the group is neither a $\pi$-group nor a $\pi'$-group. In \cite[Lemma 3.1]{Gi}, such  a simple group is characterized to be  $G=PSL(2,q)$, where $q>3$, $q(q-1)\equiv 0(3)$, $q\equiv-1(4)$ and $\pi(q+1)\subseteq \pi$, $\pi(q(q-1)/2)\subseteq \pi'$. Hence $|\pi(G)\cap \pi'|\ge |\pi(q(q-1)/2)|\ge 2$, which is not possible and proves that $G$ is $\pi$-separable.\qed

\section{General framework: The theory of soluble groups}

Before focusing on the universe of $\pi$-separable groups, and proving our main results in Section~\ref{sCarter-like}, we present briefly here the general framework on the theory of soluble groups and classes of groups, where they are relevant to our concerns. It makes clear  the origin of main concepts in our extension to $\pi$-separable groups, $\pi$ a set of primes, particularly the one of $\mathfrak N^\pi$-Dnormal subgroups, associated to the class $\mathfrak N^\pi$ of groups which are the direct product of a $\pi$-group and a nilpotent $\pi'$-group.\smallskip

A Carter subgroup of a group is a self-normalizing nilpotent subgroup. The well-known result of R.~W.~Carter states that each finite soluble group possesses exactly one conjugacy class of Carter subgroups, namely the $\mathfrak  N$-projectors (or $\mathfrak  N$-covering subgroups), for the class $\mathfrak  N$ of finite nilpotent groups.
\smallskip

We recall that given a class of groups $\mathfrak  X$, a subgroup $U$ of a group $G$ is called an $\mathfrak  X$-projector of $G$ if $UK/K$ is an $\mathfrak  X$-maximal subgroup of $G/K$ (i.e. maximal as subgroup of $G/K$ in $\mathfrak  X$)  for all $K\unlhd G$. The (possibly empty) set of $\mathfrak  X$-projectors of $G$ will be denoted by $\proj_{\mathfrak  X}(G)$.

Also, an $\mathfrak  X$-covering subgroup of $G$ is a subgroup $E$ of $G$ with the property that $E\in \proj_{\mathfrak  X}(H)$ whenever $E\le H\le G$. The set of $\mathfrak  X$-covering subgroups of $G$ will be denoted by $\cov_{\mathfrak  X}(G)$.
\smallskip

Schunck classes  $\mathfrak  X$ appear as the classes of groups for which all finite groups have associated projectors.  If the group $G$ is soluble, then $\cov_{\mathfrak  X}(G)=\proj_{\mathfrak  X}(G)$  and it is a non-empty conjugacy class of subgroups of $G$. Saturated formations are classical and relevant special cases of Schunck classes.

 Let us recall also that a class of groups $\mathfrak  X$ is a formation if every epimorphic image of a group in $\mathfrak  X$ belongs to $\mathfrak  X$, and $G/(N_1\cap N_2)\in \mathfrak  X$ whenever $N_1,N_2\unlhd G$ with $G/N_1, G/N_2\in \mathfrak  X$. In this case,  the $\mathfrak  X$-residual of a group $G$, denoted $G^{\mathfrak  X}$, is the smallest normal subgroup of $G$ with quotient group in $\mathfrak  X$; for the class $\mathfrak  X=\mathfrak  E_\pi$ of all $\pi$-groups, $\pi$ a set of primes,  $G^{\mathfrak  E_\pi}=O^\pi(G)$ is the $\pi$-residual of $G$, also described as the subgroup generated by all $\pi'$-subgroups of $G$. A formation $\mathfrak  X$ is said to be saturated if $G\in \mathfrak X$ whenever $G/\Phi(G)\in \mathfrak X$, where $\Phi(G)$ denotes the Frattini subgroup of $G$.

 It is well-known that non-empty saturated formations are exactly local formations, i.e. classes of groups $LF(f)$ consisting of the groups $G$ such that,  for every prime $p\in \p$, it holds that $G/C_G(H/K)\in f(p)$ whenever $H/K$ is a chief factor of $G$ whose order is divisible by $p$, and where $f$ is a formation function, i.e. $f(p)$ is a (possibly empty) formation for each prime $p$. The characteristic of the saturated formation is $\Char (LF(f))=\{p\in \p\mid f(p)\neq \emptyset\}$, and $LF(f)$ is said to be locally defined by the formation function $f$. We shall refer also to the fact that a local formation $\mathfrak  F$ always has a smallest local definition, i.e. a formation function $\underline{f}$ such that $\mathfrak  F=LF(\underline{f})$, and $\underline{f}(p)\subseteq g(p)$ for every prime $p$ and any other formation function $g$ such that $\mathfrak  F=LF(\underline{f})=LF(g)$.

We refer to \cite[Chapters III, IV]{DH} for details related to the theories of Schunck classes and formations.
\smallskip

Let $\pi$ be a set of primes. Let $$\mathfrak N^\pi=\mathfrak  E_\pi \times \mathfrak  N_{\pi'}=(G=H\times K\mid H\in \mathfrak  E_\pi,\ K\in \mathfrak  N_{\pi'}),$$ $\mathfrak  E_\pi$  the class of all $\pi$-groups and $\mathfrak  N_{\pi'}$ the class of all nilpotent $\pi'$-groups.

In the particular cases when either $\pi=\emptyset$ or $\pi=\{p\}$, $p$ a prime, ($|\pi|\le 1$), then $\mathfrak N^\pi=\mathfrak  N$ is the class of all nilpotent groups.
\smallskip

Our main results Theorems~\ref{teo1}~and~\ref{teo2} extend the existence and properties of Carter subgroups in soluble groups to $\rho$-separable groups, $\rho$ a set of primes, with appropriate class $\mathfrak  N^\rho$ or $\mathfrak  N^{\rho'}$ playing the role of the class $\mathfrak  N$ of nilpotent groups.
\smallskip

We shall appeal also to the concept of  $\mathfrak N^\pi$-Dnormal subgroup, as $\mathfrak N^\pi$ is a saturated formation.
\smallskip

The concept of $\mathfrak  G$-Dnormal subgroups for a non-empty saturated formation $\mathfrak  G$, which was given by K. Doerk in the universe of finite soluble groups, and appears for the first time in \cite[Definition 3.1]{AP}, is also available for arbitrary finite groups, as defined next. For notation, if $G$ is a group, $\rho$ a set of primes, $G_\rho \in \Hall_\rho(G)$ and $H\le G$, we write $G_\rho\searrow H$ to mean that $G_\rho$ reduces into $H$, i.e. $G_\rho\cap H\in \Hall_\rho (H)$.

\begin{de}\textup{\cite[Definition 3.1]{AP}} Let $\mathfrak  G$ be a non-empty saturated formation and let $G$ be a group. A subgroup $H$ of $G$ is said to be $\mathfrak  G$-Dnormal in $G$ if $\pi(|G:H|)\subseteq \Char (\mathfrak  G)$, and for every $p\in \Char (\mathfrak  G)$ it holds that
$$[H_G^p,H^{\underline{g}(p)}]\le H,$$
where $\underline{g}$ denotes the smallest local definition of $\mathfrak  G$ as local formation, and $H_G^p=\langle G_p\in \Syl_p(G)\mid G_p\searrow H\rangle$.
\end{de}

The saturated formation $\mathfrak N^\pi=\mathfrak  E_\pi \times \mathfrak  N_{\pi'}=LF(f)=LF(\underline{f})$  is locally defined by the formation function $f$ given by  $f(p)=\mathfrak  E_p$ if $p\in \pi'$, and $f(p)=\mathfrak  E_\pi$ if $p\in \pi$. Then the smallest local definition is given by $\underline{f}(p)=(1)$ if $p\in \pi'$, and $\underline{f}(p)= \begin{cases}
(1)& \text{if }\pi=\{p\},\\
\mathfrak  E_\pi & \text{if }p\in \pi,\ |\pi|\ge 2.
\end{cases}$
\medskip

Hence, for $\mathfrak N^\pi=\mathfrak  N$, $\mathfrak  N$-Dnormal subgroups are exactly normal subgroups.
\smallskip

In the case $|\pi|\ge 2$, a subgroup $H$ of a group $G$ is $\mathfrak N^\pi$-Dnormal if it satisfies the following conditions:
\begin{description}
\item[$(1)$] whenever $p\in \pi'$ and $G_p\in \Syl_p(G)$, $G_p\searrow H$, then $G_p\le N_G(H)$;
\item[$(2)$] whenever $p\in \pi$ and $G_p\in \Syl_p(G)$, $G_p\searrow H$, then $G_p\le N_G(O^\pi(H))$.

\end{description}

Note that normal subgroups are $\mathfrak N^\pi$-Dnormal for any set of primes $\pi$.

\begin{rem} Regarding the previous statement, note that for any $X\le G$ it holds that $[X,O^\pi(H)]\le H$ if and only if $X\le N_G(O^\pi(H))$.
\smallskip

{\it Proof.} Assume that $[X,O^\pi(H)]\le H$. We consider $O^\pi(H)=\langle H_q\mid H_q\in \Syl_q(H),\ q\in \pi'\rangle$. Then
$[X,O^\pi(H)]\le \langle O^\pi(H)^x\mid x\in X\rangle=\langle H_q^x\mid x\in X,\ H_q\in \Syl_q(H),\ q\in \pi'\rangle=O^\pi(H).$ The converse is clear.\qed
\end{rem}

The next proposition provides a useful characterization of $\mathfrak N^\pi$-Dnormal subgroups.
\begin{pro}\label{pro1} Let $H$ be a subgroup  of a group $G$. Then:
\begin{enumerate}\item Assume that $|\pi|\le 1$. Then $\mathfrak N^\pi =\mathfrak  N$ and $H$ is $\mathfrak  N$-Dnormal in $G$ if and only $H$ is normal in $G$.
\item Assume that $|\pi|\ge 2$. Then the following statements are equivalent:
\begin{itemize} \item[(i)] $H$ is $\mathfrak N^\pi$-Dnormal in $G$;
\item[(ii)] $O^{\pi}(H)\unlhd G$ and $O^{\pi}(G)\le N_G(H)$.
\end{itemize}
\end{enumerate}
\end{pro}
   \smallskip

{\it Proof.} Part~1 is clear. For Part~2,  since $O^\pi(G)=\langle G_p\mid G_p\in \Syl_p(G),\ p\in \pi'\rangle$, it is clear that (ii) implies (i). Conversely, assume the (i) holds, i.e. $H$ is $\mathfrak N^\pi$-Dnormal in $G$. By Sylow's theorem, for each prime $p$, there exists $G_p\in \Syl_p(G)$ such that $G_p\searrow H$. The definition of $\mathfrak N^\pi$-Dnormality implies that $G_p\le N_G(O^{\pi}(H))$ if $p\in \pi$, and $G_p\le N_G(H)\le N_G(O^{\pi}(H))$ if $p\in \pi'$. Consequently,
$G=\langle G_p\mid p\in \p\rangle\le N_G(O^{\pi}(H))$, i.e. $O^{\pi}(H)\unlhd G$. In particular, for any $p\in \pi'$ and any $G_p\in \Syl_p(G)$, it holds that  $G_p\cap H=G_p\cap O^{\pi}(H)\in \Syl_p(O^{\pi}(H))=\Syl_p(H)$, which means that $G_p\searrow H$, and then $G_p\le N_G(H)$, because $H$ is $\mathfrak N^\pi$-Dnormal in $G$. Hence, $O^\pi(G)=\langle G_p\mid G_p\in \Syl_p(G),\ p\in \pi'\rangle\le N_G(H)$, and we are done.\qed
\medskip

For notation, whenever a group $X\in \mathfrak N^\pi$, we write $X=X_{\pi}\times X_{\pi'}$ where $X_{\pi}=O_{\pi}(X)\in \mathfrak  E_{\pi}$ and $X_{\pi'}=O_{\pi'}(X)\in \mathfrak  N_{\pi'}$.
\smallskip

\begin{cor}\label{cor1} Assume that $|\pi|\ge 2$ and let $H$ be a subgroup of a group $G$ such that $H=H_\pi\times H_{\pi'}\in \mathfrak N^\pi$.
Then $H$ is $\mathfrak N^\pi$-Dnormal in $G$ if and only if $H_{\pi'}\unlhd G$ and $O^{\pi}(G)\le N_G(H)$.
\end{cor}
\smallskip

{\it Proof.} This is a consequence of Proposition~\ref{pro1}(2) since in this case $H_{\pi'}=O^{\pi}(H)$.\qed

\section{Carter-like subgroups in $\pi$-separable groups}\label{sCarter-like}

Let $\pi$ be a set of primes.  As above, set $\mathfrak N^\pi=\mathfrak  E_\pi \times \mathfrak  N_{\pi'}$, where $\mathfrak  E_\pi$ is the class of $\pi$-groups and $\mathfrak  N_{\pi'}$ is  the class of nilpotent $\pi'$-groups.
 \smallskip

 We prove in this section that if $G$ is a $\pi'$-soluble group, then $\mathfrak N^\pi$-projectors coincide with $\mathfrak N^\pi$-covering subgroups, and they form a conjugacy class of self-$\mathfrak N^\pi$-Dnormalizing subgroups of $G$ (see Definition~\ref{defselfnormal}, and  Theorems~\ref{teo1}, \ref{teo2}). If $\mathfrak N^\pi=\mathfrak  N$ is the class of nilpotent groups and $G$ is a soluble group, these are the Carter subgroups.
 \smallskip

We notice that Burnside's $p^aq^b$-theorem together with the Feit-Thompson theorem imply that $\pi$-separable groups are $\pi'$-soluble whenever  $|\pi'|\le 2$, or $2\in \pi$ if $|\pi'|\ge 3$.
\smallskip

Notice also that, by the Feit-Thompson theorem, for any set of primes $\rho$, a $\rho$-separable group is either $\rho$-soluble or $\rho'$-soluble.
\smallskip

For our main results we quote a series of results from \cite[III. Section 3]{DH} and adhere to the notation there, though specialized  to our saturated formation $\mathfrak N^\pi$ and our purposes.  We notice that  $b(\mathfrak N^\pi)\subseteq \mathcal P_1\cup\mathcal P_2$, where $b(\mathfrak N^\pi)$ is the class of groups $G\notin \mathfrak N^\pi$ but whose proper epimorphic images belong to $\mathfrak N^\pi$, and for each $i=1,2$, $\mathcal P_i$ is the class of primitive groups with a unique minimal normal subgroup, which is abelian for $i=1$, and non-abelian for $i=2$.

\begin{lem}\label{lem2} \begin{enumerate} \item \emph{\cite[III. Proposition (3.7)]{DH}} For a group $G$, whenever $N\unlhd G$, $N\le V\le G$, $U\in \proj_{\mathfrak N^\pi}(V)$, and $V/N\in \proj_{\mathfrak N^\pi}(G/N)$, then $U\in \proj_{\mathfrak N^\pi}(G)$.
\item \emph{\cite[III. Lemma (3.9)]{DH}} Assume that $G\in b(\mathfrak N^\pi)$. Then:
\begin{enumerate}\item If $G\in \mathcal P_1$, then $\cov_{\mathfrak N^\pi}(G)$ and $\proj_{\mathfrak N^\pi}(G)$ both coincide with the non-empty set comprising those subgroups of $G$ which are complements in $G$ to the minimal normal subgroup of $G$.
\item If $G\in \mathcal P_2$, then $\proj_{\mathfrak N^\pi}(G)$ is non-empty and consists of all $\mathfrak N^\pi$-maximal subgroups of $G$ which supplement the minimal normal subgroup of $G$ in $G$.
\end{enumerate}

\item \emph{\cite[III. Theorem (3.10)]{DH}} For any group $G$, $\cov_{\mathfrak N^\pi}(G)\subseteq \proj_{\mathfrak N^\pi}(G)\neq \emptyset$.

\item \emph{\cite[III. Theorem (3.14)]{DH}} Let $N$ be a nilpotent normal subgroup of a group $G$, and let $H$ be an $\mathfrak N^\pi$-maximal subgroup of $G$ such that $G=HN$. Then $H\in \proj_{\mathfrak N^\pi}(G)$.

\item \emph{(\cite[III. Theorem (3.19)]{DH},\cite{F})} Let $\mathfrak B$ be the class of all $\pi'$-soluble groups. The statement ``\,$\proj_{\mathfrak N^\pi}(G)$ is a conjugacy class of $G$'' is true for all groups $G\in \mathfrak B$ if and only if it is true for all $G\in b(\mathfrak N^\pi)\cap \mathfrak B$.
    \item \emph{(\cite[III. Remark (3.20)(b)]{DH},\cite{F})} Let $\mathfrak B$ be the class of all $\pi'$-soluble groups. The statement ``\,$\proj_{\mathfrak N^\pi}(G)=\cov_{\mathfrak N^\pi}(G)$'' is true for all groups $G\in \mathfrak B$ if and only if it is true for all $G\in b(\mathfrak N^\pi)\cap \mathfrak B$.
\end{enumerate}
\end{lem}

We still quote the following result for our purposes.

\begin{lem}\label{lem2+}\emph{\cite[Theorems 4.1.18, 4.2.17]{BE}} Let $\mathfrak  H$ be a saturated formation and let $G$ be a group  whose $\mathfrak  H$-residual $G^{\mathfrak  H}$ is abelian. Then $G^{\mathfrak  H}$ is complemented in $G$, any two complements are conjugate in $G$, and the complements are the $\mathfrak  H$-projectors of $G$.
\end{lem}

\begin{lem}\label{lem3} Let $M=M_{\pi}\times M_{\pi'}$ be an $\mathfrak N^\pi$-maximal subgroup of a $\pi$-separable group $G$. Then:
\begin{enumerate}\item $M=M_{\pi'}C_G(M_{\pi'})_{\pi}$ for some $C_G(M_{\pi'})_{\pi}\in Hall_{\pi} (C_G(M_{\pi'}))$.
\item If $H=H_{\pi}\times H_{\pi'}$ is another $\mathfrak N^\pi$-maximal subgroup of $G$ and $M_{\pi'}^x= H_{\pi'}$ for some $x\in G$, then $M^g=H$ for some $g\in G$.
\end{enumerate}
\end{lem}

\smallskip

{\it Proof.} 1. We have that $M_{\pi}\le C_G(M_{\pi'})$ and so $M_{\pi}\le C_G(M_{\pi'})_{\pi}$ for some $C_G(M_{\pi'})_{\pi}\in Hall_{\pi} (C_G(M_{\pi'}))$. But $M_{\pi'}C_G(M_{\pi'})_{\pi}\in \mathfrak N^\pi$, which implies that  $M=M_{\pi'}C_G(M_{\pi'})_{\pi}\in \mathfrak N^\pi$ by the maximality of $M$.
\smallskip

2. The hypothesis implies that $C_G(M_{\pi'})^x= C_G(M_{\pi'}^x) =
C_G(H_{\pi'}) $. Then $M_{\pi}^x, H_{\pi}\in \Hall_{\pi}(C_G(H_{\pi'}))$ and $M_{\pi}^{xy}=H_{\pi}$ for some $y\in C_G(H_{\pi'})$. Consequently,
$M^{xy}= M_{\pi}^{xy}M_{\pi'}^{xy}=H_{\pi}H_{\pi'}=H$, and we are done.\qed

\begin{teo}\label{teo1} If $G$ is a $\pi'$-soluble group, then $\emptyset \neq \proj_{\mathfrak N^\pi}(G)= \cov_{\mathfrak N^\pi}(G)$ and it is a conjugacy class of $G$.
\end{teo}
\smallskip

{\it Proof.} By Lemma~\ref{lem2}, parts~(3),~(5),~(6), we may assume that $G\in b(\mathfrak N^\pi)$. Since $\mathfrak N^\pi$ is a saturated formation, $G\in \mathcal P_1\cup \mathcal P_2$; let $N$ be the minimal normal subgroup of $G$.

If $G\in \mathcal P_1$, the result follows by Lemmas~\ref{lem2}(2)(a)~and~\ref{lem2+}.

Assume now that $G\in \mathcal P_2$.  We know by Lemma~\ref{lem2}(2)(b) that $\proj_{\mathfrak N^\pi}(G)$ is non-empty and consists of all $\mathfrak N^\pi$-maximal subgroups of $G$ which supplement $N$ in $G$. We prove first that these subgroups are conjugate in $G$.

Let $M=M_{\pi}\times M_{\pi'}$, $H=H_{\pi}\times H_{\pi'}$ be $\mathfrak N^\pi$-maximal subgroups of $G$ such that $G=NM=NH$. Since $G$ is $\pi'$-soluble,  and $N$ is non-abelian, $N$ is a $\pi$-group and, consequently, $M_{\pi'}, H_{\pi'}\in \Hall_{\pi'}(G)$. Hence, there exists $x\in G$ such that $M_{\pi'}^x=H_{\pi'}$ and $M$ and $H$ are conjugate by Lemma~\ref{lem3}(2).
\smallskip

We claim now that $\proj_{\mathfrak N^\pi}(G)\subseteq \cov_{\mathfrak N^\pi}(G)$, which will conclude the proof.

Let $M\in \proj_{\mathfrak N^\pi}(G)$, i.e. $M=M_{\pi}\times M_{\pi'}$ is an $\mathfrak N^\pi$-maximal subgroup of $G$ such that $G=NM$. Let $M\le L\le G$. We aim to prove that $M\in \proj_{\mathfrak N^\pi}(L)$. Let $T\in \proj_{\mathfrak N^\pi}(L)$. We notice that $L=M(L\cap N)$. Then $L/(L\cap N)\cong M/(M\cap N)\in \mathfrak N^\pi$. Since $T(L\cap N)/(L\cap N)$ is  $\mathfrak N^\pi$-maximal in $L/(L\cap N)$, it follows that $L=T(L\cap N)$. Hence $T_{\pi'}, M_{\pi'}\in \Hall_{\pi'}(L)$ and, moreover, $T$ and $M$ are $\mathfrak N^\pi$-maximal subgroups of $L$. By Lemma~\ref{lem3}(2), $T$ and $M$ are conjugate in $L$ and $M\in \proj_{\mathfrak N^\pi}(L)$.\qed

\begin{rem} In Theorem~\ref{teo1} the hypothesis of $\pi'$-solubility cannot be weakened to $\pi$-separability.  Otherwise, for the particular case when $\pi=\emptyset$, the result would hold for every finite group and the formation $\mathfrak N^\pi=\mathfrak  N$ of nilpotent groups, which is not true. Particularly, also if $\pi\neq \emptyset$, one can consider for instance $\pi=\p - \{2,3,5\}$, $\pi'=\{2,3,5\}$ and $G=\Alt(5)$ the alternating group of degree $5$. The group $G$ is obviously $\pi$-separable, the $\mathfrak  N^{\pi}$-projectors are the $\mathfrak  N$-projectors, which do not form a conjugacy class of subgroups, as they are all the Sylow subgroups of $G$; and $G$ has no $\mathfrak  N$-covering subgroups.
\end{rem}

\begin{de}\label{defselfnormal}
A subgroup $H$ of a group $G$ is said to be \emph{self-$\mathfrak N^\pi$-Dnormalizing} in $G$ if whenever $H\le K\le G$ and $H$ is $\mathfrak N^\pi$-Dnormal in $K$, then $H=K$.
\end{de}

We prove next that $\mathfrak N^\pi$-projectors are self-$\mathfrak N^\pi$-Dnormalizing subgroups.

\begin{pro}\label{prop2} Let $H$ be an $\mathfrak N^\pi$-projector of a  $\pi'$-soluble group $G$. Then $H$ is self-$\mathfrak N^\pi$-Dnormalizing in $G$.
\end{pro}
\smallskip

{\it Proof.} Assume that $H\le K\le G$ and $H$ is $\mathfrak N^\pi$-Dnormal in $K$. We aim to prove that $H=K$. By Theorem~\ref{teo1} and Corollary~\ref{cor1}, we have that $H\in \cov_{\mathfrak N^\pi}(K)$ and $H_{\pi'}\unlhd K$. Then $H/H_{\pi'}\in \proj_{\mathfrak N^\pi}(K/H_{\pi'})$   and $H/H_{\pi'}\le K_{\pi}H_{\pi'}/H_{\pi'}\in \mathfrak N^\pi$, for any $K_{\pi}\in \Hall_{\pi}(K)$, which implies that $H=K_{\pi}H_{\pi'}$. Whence, if $K_{\pi'}\in \Hall_{\pi'}(K)$, then $K=HK_{\pi'}$ and $H\unlhd K$ by Proposition~\ref{pro1}. If $H< K$, then $H< HK_p$ for some $p\in \pi'$ and $1\ne K_p\in \Syl_p(K)$. But $HK_p/H\in \mathfrak N^\pi$ which contradicts the fact that $H\in \cov_{\mathfrak N^\pi}(K)$.\qed

\begin{lem}\label{lem4} Assume that $|\pi|\ge 2$, $H=H_{\pi}\times H_{\pi'},K=K_{\pi}\times K_{\pi'}\in \mathfrak N^\pi$ and $H\le K$. Then $H$ is $\mathfrak N^\pi$-Dnormal in $K$ if and only if $H_{\pi'}\unlhd K_{\pi'}$.
\end{lem}
\smallskip

{\it Proof.} We notice that $[K_{\pi},H_{\pi'}] =1$ and $[K_{\pi'},H_{\pi}] =1$. Consequently, Corollary~\ref{cor1} implies that $H$ is $\mathfrak N^\pi$-Dnormal in $K$ if and only if $H_{\pi'}\unlhd K_{\pi'}$.\qed

\begin{pro} Assume that $H=H_{\pi}\times H_{\pi'}< L=L_{\pi}\times L_{\pi'}\in \mathfrak N^\pi$. Then there exists $K\le L$ such that $H<K$ and $H$ is $\mathfrak N^\pi$-Dnormal in $K$.
\end{pro}
\smallskip

{\it Proof.} If $|\pi|\le 1$, then $\mathfrak N^\pi = \mathfrak  N$, and the result is clear. In the case $|\pi|\ge 2$, by Lemma~\ref{lem4}, if $H_{\pi'}=L_{\pi'}$, then $H$ is $\mathfrak N^\pi$-Dnormal in $L$, and we are done. Otherwise, there exists $T\le L_{\pi'}$ such that $H_{\pi'}\vartriangleleft T$, i.e. $H_{\pi'}$ is a proper normal subgroup of $T$, because $L_{\pi'}$ is nilpotent. We can consider now the subgroup $K=L_{\pi}T\le L$ which satisfies that $H<K$ and $H$ is $\mathfrak N^\pi$-Dnormal in $K$, which concludes the proof.\qed
\smallskip

As a consequence we can state the following:

\begin{cor}\label{cor2} If $H\in \mathfrak N^\pi$ is a self-$\mathfrak N^\pi$-Dnormalizing subgroup of a group $G$, then $H$ is $\mathfrak N^\pi$-maximal in $G$.
\end{cor}
\begin{rem} It is not true in general that $\mathfrak N^\pi$-projectors  of $\pi'$-soluble groups are exactly self-$\mathfrak N^\pi$-Dnormalizing subgroups in $\mathfrak N^\pi$. Otherwise,  $\mathfrak N^\pi\cap \mathfrak  S$ would be either $\mathfrak  N$ or $\mathfrak  S$, the class of all soluble groups, by \cite[Proposition 4.1]{AP04}. But we see next that a corresponding result to \cite[Theorem 4.2]{AP04} is still possible. That reference provides a corresponding result to our next Theorem~\ref{teo2}, for finite soluble groups, subgroup-closed saturated formations and associated projectors.
\end{rem}

As a consequence of Lemma~\ref{lem2+} we can state the following.

\begin{lem}\label{lem5} Let $\mathfrak  H$ be a saturated formation, $X$ be a group and $H$ be an $\mathfrak  H$-projector of $X$. Then $H\cap X^{\mathfrak  H}\le (X^{\mathfrak  H})'$.
\end{lem}

\begin{teo}\label{teo2}  For a subgroup $H$ of a $\pi'$-soluble group $G$ the following statements are pairwise equivalent:
\begin{enumerate}\item $H$ is an $\mathfrak N^\pi$-projector of $G$.
\item $H$ is an $\mathfrak N^\pi$-covering subgroup of $G$.
\item $H\in \mathfrak N^\pi$ is a self-$\mathfrak N^\pi$-Dnormalizing subgroup of $G$ and $H$ satisfies the following property:
\begin{gather}\text{If}\  H\le X\le G,\ \text{then}\ H\cap X^{\mathfrak N^\pi}\le (X^{\mathfrak N^\pi})'.\tag{*}
\end{gather}

\end{enumerate}
\end{teo}
\smallskip

{\it Proof.} The equivalence $1 \leftrightarrow 2$ has been proven in Theorem~\ref{teo1}. On the other hand, Proposition~\ref{prop2} and Lemma~\ref{lem5} prove $2\rightarrow 3$. We prove next that $3\rightarrow 1$.
\smallskip

Let  $H\in \mathfrak N^\pi$ be a self-$\mathfrak N^\pi$-Dnormalizing subgroup of $G$ satisfying property $(^*)$. We aim to prove that $H\in \proj_{\mathfrak N^\pi}(G)$. We notice that $H$ is $\mathfrak N^\pi$-maximal in $G$ by Corollary~\ref{cor2}.
\smallskip

If $G\in \mathfrak N^\pi$, then $H=G$ and the result follows. So that we may assume that $G\notin \mathfrak N^\pi$. We argue by induction on the order of $G$. Let $N$ be a minimal normal subgroup of $G$ such that $N\le G^{\mathfrak N^\pi}$.
\smallskip

We distinguish the following cases:
\begin{description}
\item[Case 1.]  $G=HN$.
\item[Case 2.] $HN<G$.
\end{description}
\smallskip

\noindent
{\bf Case 1.}  If $N$ is abelian, the result follows by Lemma~\ref{lem2}(4). Assume that $N$ is not abelian. Let $K\unlhd G$ such that $G/K\in b(\mathfrak N^\pi)$. Then $N$ is not contained in $K$ because $G/N\cong H/(H\cap N)\in \mathfrak N^\pi$. In particular, $N\cap K=1$ and $N\cong NK/K$ is a minimal normal subgroup of $G/K=(NK/K)(HK/K)\in \mathcal P_2$, with $HK/K<G/K$ because $HK/K\in \mathfrak N^\pi$. By Lemma~\ref{lem2}(2)(b), $HK/K\le P/K$ for some $P/K\in \proj_{\mathfrak N^\pi}(G/K)$.
We have now that $H\le P<G$. The inductive hypothesis implies that $H\in \proj_{\mathfrak N^\pi}(P)$, and from Lemma~\ref{lem2}(1), $H\in \proj_{\mathfrak N^\pi}(G)$, as claimed.
\smallskip

\noindent
{\bf Case 2.}  In this case $HN<G$ and the inductive hypothesis implies that $H\in \proj_{\mathfrak N^\pi}(HN)$. We prove first that $HN/N$ satisfies property $(^*)$ in $G/N$. Assume that $HN/N\le X/N\le G/N$. If $X<G$, then $H\in \proj_{\mathfrak N^\pi}(X)$ by inductive hypothesis and then $HN/N\cap (X/N)^{\mathfrak N^\pi}\le ((X/N)^{\mathfrak N^\pi})'$.
Otherwise, $X=G$ and so $(HN/N)\cap (X/N)^{\mathfrak N^\pi}= (HN/N)\cap (G/N)^{\mathfrak N^\pi}= (HN/N)\cap G^{\mathfrak N^\pi}/N= (H\cap G^{\mathfrak N^\pi})N/N\le (G^{\mathfrak N^\pi})'N/N= ((G/N)^{\mathfrak N^\pi})'=((X/N)^{\mathfrak N^\pi})'.$

We claim that $HN/N\in \mathfrak N^\pi$ is self-$\mathfrak N^\pi$-Dnormalizing in $G/N$. Then the result follows by inductive hypothesis together with Lemma~\ref{lem2}(1).
\smallskip

Assume that $HN/N$ is $\mathfrak N^\pi$-Dnormal in $L/N\le G/N$.

If $L<G$, the inductive hypothesis implies that $H$ is an $\mathfrak N^\pi$-projector of $L$ and so $HN=L$ by Proposition~\ref{prop2}. So that we may assume that $L=G$ and $HN/N$ is $\mathfrak N^\pi$-Dnormal in $ G/N$.
\smallskip

We split the rest of the proof into the following steps:
\smallskip

\noindent
{\it Step 1.} $N_G(HN)=HN$.
\smallskip

If $g\in N_G(HN)$, then $HN=H^gN$. Since $H\in \proj_{\mathfrak N^\pi}(HN)$, and this is a conjugacy class of subgroups of $HN$ by Theorem~\ref{teo1}, $H^g=H^x$ for some $x\in HN$. Consequently, $gx^{-1}\in N_G(H)=H$, and $g\in HN$.
\smallskip

\noindent
{\it Step 2.} $O^{\pi}(G)\le HN$.
\smallskip

Since $HN/N$ is $\mathfrak N^\pi$-Dnormal in $G/N$, we have by Proposition~\ref{pro1} and Step~1 that $O^{\pi}(G)N/N\le N_{G/N}(HN/N)=HN/N$.
\smallskip

\noindent
{\it Step 3.} $O^{\pi}(G)H=HN$ and it is a maximal subgroup of $G$.
\smallskip

Assume that $O^{\pi}(G)\le HN\le T< G$. The inductive hypothesis implies that $H\in \proj_{\mathfrak N^\pi}(T)$. Moreover, $T/O^{\pi}(G)\in \mathfrak  E_{\pi}\subseteq \mathfrak N^\pi$. Hence $T=O^{\pi}(G)H=HN$.
\smallskip

\noindent
{\it Step 4.} $G/N= (H_{\pi'}N/N)\times O_{\pi}(G/N)\in \mathfrak N^\pi$. In particular, $G^{\mathfrak N^\pi}=N$.
\smallskip

We have that $H_{\pi'}N/N\unlhd G/N$ and $H_{\pi'}N/N\in \Hall_{\pi'}(G/N)$ by Corollary~\ref{cor1} and Step~2. Since $[H_{\pi},H_{\pi'}]=1$ and $G/N$ is $\pi$-separable, and equivalently $\pi'$-separable, it follows  by \cite[6. Theorem 3.2]{Go} that \begin{gather*}(H_\pi N/N)O_\pi (G/N)/O_\pi (G/N)\le C_{(G/N)/O_\pi(G/N)}(O_{\pi'}((G/N)/O_\pi(G/N)))\\ \le O_{\pi'}((G/N)/O_\pi(G/N)),\end{gather*} which implies  $H_{\pi}N/N\le O_{\pi}(G/N)$.

If $H_{\pi}N/N= O_{\pi}(G/N)$, then $HN/N\unlhd G/N$ and $G=HN$ by Step~1, a contradiction.

Consequently, we may assume that $H_{\pi}N/N< O_{\pi}(G/N)$. Then $HN/N< (H_{\pi'}N/N) O_{\pi}(G/N)\le G/N$. By Step~3, $(H_{\pi'}N/N) O_{\pi}(G/N)= G/N \in \mathfrak N^\pi$.
\smallskip

\noindent
{\it Step 5.}  $G=NP$ where $P=P_\pi\times P_{\pi'}\in \proj_{\mathfrak N^\pi}(G)$.
\smallskip

It follows by Step~4.
\smallskip

\noindent
{\it Step 6.}  $N\in \mathfrak  E_{\pi'}$.
\smallskip

If $N\in \mathfrak  E_{\pi}$, then $H_{\pi'}, P_{\pi'}\in \Hall_{\pi'}(G)$. Since $H$ and $P$ are both $\mathfrak N^\pi$-maximal subgroups, it follows by Lemma~\ref{lem3} that $H=P^x\in \proj_{\mathfrak N^\pi}(G)$ and then $HN=G$, a contradiction. Since $G$ is $\pi$-separable, $N\in \mathfrak  E_{\pi'}$.
\smallskip

\noindent
{\it Step 7.}  Final contradiction.
\smallskip

Steps~4,~6, imply that $N=G^{\mathfrak N^\pi}$ is abelian, since $G$ is $\pi'$-soluble. By hypothesis, $H\cap N=1$. Moreover, $P_{\pi'}N=H_{\pi'}N= O_{\pi'}(G)$ and $P_\pi\in \Hall_{\pi}(G)$. There is no loss of generality to assume that $H_\pi\le P_{\pi}$ and then $[H_\pi, P_{\pi'}]=1$.

Since $H$ is an $\mathfrak N^\pi$-maximal subgroup of $HN=H_{\pi}H_{\pi'}N$, we deduce that $H_{\pi'}$ is an $\mathfrak  N$-maximal subgroup of $C_{H_{\pi'}N}(H_{\pi})=H_{\pi'}C_N(H_{\pi})$. But $C_N(H_{\pi})\unlhd C_{H_{\pi'}N}(H_{\pi})$ and $C_N(H_{\pi})$ is nilpotent. By Lemma~\ref{lem2}(4), $H_{\pi'}\in \proj_{\mathfrak  N}(C_{H_{\pi'}N}(H_{\pi}))$. On the other hand, $$C_{H_{\pi'}N}(H_{\pi})=C_{P_{\pi'}N}(H_{\pi})=P_{\pi'}C_{N}(H_{\pi}).$$ Since $P_{\pi'}$ is nilpotent, it follows again by Lemma~\ref{lem2}(4), that $P_{\pi'}\le H_{\pi'}^x\in \proj_{\mathfrak  N}(C_{H_{\pi'}N}(H_{\pi}))$, for some $x\in C_N(H_\pi)$. If
$P_{\pi'}< H_{\pi'}^x$, since $P_{\pi'}N=H_{\pi'}^xN$, we have that $H_{\pi'}^x\cap N\neq 1$ and so also $H\cap N\neq 1$, a contradiction. Therefore,
$P_{\pi'}= H_{\pi'}^x$, and Lemma~\ref{lem3} implies that $H^g=P$ for some $g\in G$ as $H$ and $P$ are $\mathfrak N^\pi$-maximal subgroups of $G$. It follows that $H\in \proj_{\mathfrak N^\pi}(G)$ and $HN=G$, the final contradiction.\qed

\bigskip

\noindent
{\bf Acknowledgments.} Research supported by Proyectos PROMETEO/2017/ 057 from the Generalitat Valenciana (Valencian Community, Spain), and PGC2018-096872-B-I00 from the Ministerio de Ciencia, Innovaci\'on y Universidades, Spain, and FEDER, European Union. The fourth author acknowledges with thanks the financial support of the Universitat de Val\`encia as research visitor (Programa Propi d'Ajudes a la Investigaci\'o de la Universitat de Val\`encia, Subprograma Atracci\'o de Talent de VLC-Campus, Estades d'investigadors convidats (2019)).

\bigskip

\noindent Milagros Arroyo-Jord\'a, Paz Arroyo-Jord\'a\\
Escuela T\'{e}cnica Superior de Ingenieros Industriales\\
Instituto Universitario de Matem\'{a}tica Pura y  Aplicada IUMPA\\
Universitat Polit\`ecnica de Val\`encia,\\
 Camino de Vera, s/n,  46022 Valencia, Spain\\
E-mail: marroyo@mat.upv.es, parroyo@mat.upv.es\\
 \\
Rex Dark\\
School of Mathematics, Statistics and Applied Mathematics,\\
National University of Ireland, University Road, Galway,
  Ireland\\
E-mail: rex.dark@nuigalway.ie\\
 \\
 Arnold D.~Feldman\\
Franklin and Marshall College, Lancaster, PA 17604-3003,
  U.S.A.\\
E-mail: afeldman@fandm.edu\\
 \\
Mar\'{\i}a Dolores P\'erez-Ramos\\
Departament de Matem\`{a}tiques, Universitat de Val\`{e}ncia,\\
C/ Doctor Moliner 50, 46100 Burjassot
(Val\`{e}ncia), Spain\\
E-mail: Dolores.Perez@uv.es

\end{document}